\documentclass{amsart}

\numberwithin{equation}{section}

\usepackage{amssymb}
\usepackage{enumerate, xspace,epsfig}
\usepackage{amscd}
\usepackage{hyperref}
\usepackage{cite}

\newtheorem{thmm}{Theorem}
\newtheorem{thm}{Theorem}[section]
\newtheorem{lem}[thm]{Lemma}

\newtheorem{defi}{Definition}

\newtheorem{rem}[thm]{Remark}

\newcommand\Co{{\mathcal C}}

\newcommand\F{{\mathcal F}}

\newcommand\J{{\mathcal J}}

\newcommand\M{{\mathcal M}}

\newcommand\V{{\mathcal V}}

\newcommand\NN{{\mathbb N}}
\newcommand\RR{{\mathbb R}}

\newcommand\ZZ{{\mathbb Z}}

\newcommand\eps{\epsilon}
\newcommand\vf{\varphi}
\newcommand\bsigma{{\bf \sigma}}

\newcommand{\ee}{\text{\it\bfseries e}}
\newcommand{\ii}{\text{\it\bfseries i}}
\newcommand{\jj}{\text{\it\bfseries j}}
\newcommand{\kk}{\text{\it\bfseries k}}
\newcommand{\xx}{\text{\it\bfseries x}}
\newcommand{\nn}{\text{\bfseries 0}}

\newcommand{\const}{\operatorname{const}}
\newcommand{\sign}{\operatorname{sign}}
\newcommand{\card}{\operatorname{\#}}
\newcommand{\eg}{\emph{e.g.}\xspace}
\newcommand{\ie}{\emph{i.e.}\xspace}
\newcommand{\cf}{\emph{cf.}\xspace}

\newcommand{\AC}{AC_\Omega}
\newcommand{\ACC}{ACC_\Omega}
\newcommand{\BV}{BV_\Omega}
\newcommand{\Ts}{S}
\newcommand{\var}{\operatorname{{var}}}
\newcommand{\Var}{\operatorname{{Var}}}

\makeatletter
\newcommand{\subjclassname@NEW}{2000 Mathematics Subject Classification}
\makeatother

\begin{document}
\title[Phase transitions in CML]{Phase transitions in a piecewise expanding
  coupled map lattice with linear nearest neighbour coupling} \author{Jean-Baptiste Bardet and
  Gerhard Keller}
\address{J.-B. Bardet, IRMAR/UFR Math\'ematiques, Universit\'e de Rennes 1,
Campus de Beaulieu, 35042 Rennes Cedex, France;\; G. Keller,
Mathematisches Institut, Universit\"at Erlangen-N\"urnberg,
  Bismarckstr. 1 1/2, 91054 Erlangen, Germany} \email{
  {\tt jean-baptiste.bardet@univ-rennes1.fr,
    keller@mi.uni-erlangen.de}} \date{\today} \thanks{G.K. thanks the
  colleagues at the UFR
  Math\'ematiques of the University of Rennes 1 for their hospitality during
  his stay in March and April 2006. This note would
  not have been written, if both authors had not had the chance to attend the
  workshop
``COUPLED MAP LATTICES 2004'' at the IHP, Paris. Also many discussions 
over the years with Carlangelo Liverani helped to shape ideas.}
\begin{abstract}
  We construct a mixing continuous piecewise linear map on $[-1,1]$ with the
  property that a two-dimensional lattice made of these maps with a linear
  north and east nearest neighbour coupling admits a phase transition. We also
  provide a modification of this construction where the local map is an
  expanding analytic circle map. The basic strategy is borroughed from
  \cite{GiMcK-2000}, namely we compare the dynamics of the CML to those of a
  probabilistic cellular automaton of Toom's type, see \cite{MacKay-2005} for
  a detailed discussion.
\end{abstract}
\subjclass[NEW]{37L40,37L60,82C20} \maketitle

\section{Introduction}
\label{sec:intro}
The purpose of this article
is to construct a continuous piecewise
linear map $\tau$ on $I=[-1,1]$ such that the coupled map lattice (CML)
$S_\epsilon:I^\Lambda\to I^\Lambda$
($\Lambda=\ZZ^2$ or $=\ZZ/(L\ZZ)^2$) 
defined by
\begin{equation}
  \label{eq:def-Seps}
  (S_\epsilon(\xx))_\ii=(1-\epsilon)\tau(x_\ii)+\frac\epsilon2(\tau(x_{\ii+\ee_1})+\tau(x_{\ii+\ee_2}))
\end{equation}
($\ee_1,\ee_2$ are the canonical unit vectors in $\Lambda$) has a phase
transition in the following sense: there are $0<\epsilon_1<\epsilon_2<\eta$ 
such that
\begin{itemize}
\item for $0<\epsilon\leq\epsilon_1$ the infinite system and also the
  finite ones have a unique invariant
  probability measure with absolutely continuous finite-dimensional
  conditional marginals (so the measure is absolutely
  continuous if $\Lambda$ is finite),
\item for $\epsilon_2\leq\epsilon\leq\eta$ the infinite system
  has at least two such invariant probability measures
  while the finite systems still have a unique absolutely continuous
  invariant probability measure.
\end{itemize}
Using the notations
\begin{eqnarray}
  \label{eq:def-Phieps}
   \Phi_\epsilon:\Omega\to\Omega,&\quad
  (\Phi_\epsilon(\xx))_\ii
  &=(1-\epsilon)x_\ii+\frac\epsilon2(x_{\ii+\ee_1}+x_{\ii+\ee_2})\;,\\
  T:\Omega\to\Omega,&\quad(T\xx)_\ii&=\tau(x_\ii)\;,
\end{eqnarray}
we can write $S_\epsilon=\Phi_\epsilon\circ T$, and as
$S_\epsilon\circ\Phi_\epsilon=\Phi_\epsilon\circ(T\circ\Phi_\epsilon)$, it is
equivalent to study instead of $S_\epsilon$ the system
\begin{equation}
  \label{eq:def-Teps}
  T_\epsilon:I^\Lambda\to I^\Lambda,\quad T_\epsilon(\xx)=T(\Phi_\epsilon(\xx))\;.
\end{equation}

In order to state our results we adopt the notation of \cite{KeLi-2005b}: Let
$\Omega=I^\Lambda$ and let $\M(\Omega)$ be the set of signed Borel measures on
$\Omega$.\footnote{We use the product topology on $\Omega$.}
We need to introduce the concept of \emph{measures of bounded
  variation}. For $\mu\in\M(\Omega)$ let
\begin{equation}\label{eq:variation}
\begin{split}
\Var\mu
&:=\sup_{\ii\in\Lambda}\sup_{|\vf|_{\Co^0(\Omega)}\leq 1}\mu(\partial_\ii\vf)
\end{split}
\end{equation}
Here $\partial_\ii$ denotes the partial derivative with respect to
$x_\ii$.\footnote{Here and in the sequel all test functions $\vf:\Omega\to\RR$
  depend on only finitely many coordinates and are $C^1$ with respect to these
  coordinates.} It is easy to prove that the set
$\BV:=\{\mu\in\M(\Omega)\;:\;\Var\mu<\infty\}$ consists of measures whose
finite dimensional marginals are absolutely continuous with respect to
Lebesgue and the density is a function of bounded variation \cite{KeLi-2004}.
In addition, such measures have finite entropy density with respect to
Lebesgue \cite[Corollary 4.1]{KeLi-2005a}.  In fact, ``$\Var$'' is a norm and,
with this norm, $\BV$ is a Banach space. Remark that for finite $\Lambda$ and
$d\mu=f\,dm^\Lambda$ where $m$ denotes Lebesgue measure on $I$, $\Var(\mu)$ is
just $\var(f)$ in the sense of functions of bounded variation.\footnote{See
  \cite{KeLi-2005a} for a careful discussion of bounded variation in the
  present context and the relevant associated properties.}

Our main result is

\begin{thmm}\label{theo:main}
  With the piecewise linear map $\tau$ we propose in section
  \ref{sec:local-map}, the system $(T_\epsilon,\Omega)$ has the following
  properties: There are $0<\epsilon_1<\epsilon_2<\eta$ such that the following
  hold:
  \begin{enumerate}[a)]
  \item\label{theo:main-a} For $\epsilon\in[0,\frac14]$, the map
    $T_\epsilon$ has at least one invariant probability measure in
    $\BV$ which is also translation invariant.
  \item\label{theo:main-b} For $\epsilon\in[0,\epsilon_1]$, the
    map $T_\epsilon$ has a unique invariant probability measure in
    $\BV$. (This measure is necessarily also translation
    invariant.)
  \item\label{theo:main-c} For $\epsilon\in[\epsilon_2,\eta]$, the map
    $T_\epsilon$ has at least two invariant probability measures
    $\mu_\epsilon^+$ and $\mu_\epsilon^-$ in $\BV$ with
    $\mu_\epsilon^+\{x_\nn\leq0\}=\mu_\epsilon^-\{x_\nn\geq0\}<\frac12$.
\end{enumerate}
\end{thmm}
Assertions \ref{theo:main-a}) and \ref{theo:main-b}) follow rather
directly from known results, essentially from \cite{KeKue-1992} and
\cite{KeLi-2005a}, respectively. For the proof of assertion
\ref{theo:main-c}) we rely heavily on the construction of a phase
transition in Toom's probabilistic cellular automaton (PCA) as presented in
\cite{LeMaSp-1990}.

The idea to link the dynamics of a CML to those of a PCA was
introduced by Gielis and MacKay in \cite{GiMcK-2000} where they
construct CMLs with simple piecewise linear Markov maps as local units
and with discrete couplings in such a way that their dynamics are
isomorphic to those of certain PCAs. This is definitively not the case
in our model because we use a ``traditional'' directed nearest
neighbour coupling, and so a number of additional arguments are
necessary to link our CML to a PCA of Toom's type.  These arguments
are provided in sections~\ref{sec:global-coupling}
and~\ref{sec:exp-estimate}.  On the other hand our local map $\tau$
has a very large number of monotone branches so that it acts nearly as
an instantaneous local random generator unlike the rather simple local
maps  
for which Boldrighini et al \cite{BBCSP-2001} exhibited
phase transitions numerically even in a one-dimensional lattice. (See
also earlier references given in that paper or in \cite{JuSchm-2005}.)

We fix some further notation: For $\Lambda'\subseteq\Lambda$ let
$\F_{\Lambda'}$ be the $\sigma$-algebra on $\Omega$ generated by the
coordinates $x_\ii$, $\ii\in\Lambda'$. 
\begin{defi}
  \begin{enumerate}[a)]
  \item The probability measure $\nu$ on $\Omega$ belongs to the class
    $\AC$ if it has absolutely continuous finite-dimensional marginal
    distributions, \ie if, for each finite $\Lambda'\subset\Lambda$,
    the projection of $\nu$ to $I^{\Lambda'}$ is absolutely continuous
    with respect to Lebesgue measure on $I^{\Lambda'}$.
  \item The probability measure $\nu$ on $\Omega$ belongs to the class
    $\ACC$ if, for each finite $\Lambda'\subset\Lambda$, $\nu$ has an
    $\F_{\Lambda\setminus\Lambda'}$-measurable family of conditional
    probability distributions on $I^{\Lambda'}$ which are all
    absolutely continuous with respect to Lebesgue measure on
    $I^{\Lambda'}$.
  \item By $\BV(T_\epsilon)$, $\AC(T_\epsilon)$, and
    $\ACC(T_\epsilon)$ we denote the $T_\epsilon$-invariant measures
    in $\BV$, $\AC$, and $\ACC$, respectively.
  \end{enumerate}
\end{defi}
\begin{rem}
We have the following two inclusions:
\begin{equation}
  \label{eq:ACC-inclusion}
  \{\nu\in\BV:\text{ $\nu$ is a probability measure}\}\subseteq\ACC\subseteq\AC
\end{equation}
The second one is obvious, the first one was proved in \cite[Lemma
4]{Keller-2000}. In our particular setting of one-sided directed couplings
even more is known. The invariant measures belong to $\ACC$, namely:
\begin{equation}
  \label{eq:ACCinv-inclusion}
  \{\nu\in\BV(T_\epsilon):\text{ $\nu$ is a probability measure}\}=\ACC(T_\epsilon)\subseteq\AC(T_\epsilon)\,,
\end{equation}
see \cite[Proposition 2(a)]{KeZw-2002}.  It follows that the unique invariant
measure from Theorem~\ref{theo:main}\ref{theo:main-b} is indeed unique within
the class $\ACC$.\footnote{Of course such statements need an assumption on the
  coupling strength $\epsilon$. In \cite{KeZw-2002} it is required that
  $|\epsilon|$ is ``sufficiently small''. What is needed for
  \eqref{eq:ACCinv-inclusion} is a suitable Lasota-Yorke inequality, and that
  this holds for $|\epsilon|\leq\frac14$ follows just like our
  Lemma~\ref{lem:lasota-yorke}.}  Whether this is true more generally sems to
be an open question.
\end{rem}  

Our goal is to construct a coupled map lattice with a phase transition and not
merely with a bifurcation at some critical parameter. As a criterion for a
true phase transition in this sense MacKay \cite{MacKay-2005} suggests that
the dynamical system $(T_\epsilon,\Omega)$ should satisfy a kind of space-time
specification property called \emph{indecomposability} - also after the phase
transition has occured. But since in our case the dynamics do not even have a
tractable symbolic representation, this property seems inappropriate here.
Instead we prove that each finite lattice version of the system
$(S_\epsilon,\Omega)$ has a unique absolutely continuous invariant measure
which is mixing and equivalent to the finite-dimensional Lebesgue measure on
$I^\Lambda$. More precisely:
\begin{thmm}
  \label{theo:finite-dim}
  Let $\Lambda=(\ZZ/L\ZZ)^2$ for some $L\in\ZZ_+$ and define
  $T_\epsilon:I^\Lambda\to I^\Lambda$ as in \eqref{eq:def-Teps} with periodic
  boundary conditions.  For the map $\tau$ from Theorem~\ref{theo:main} and
  $\epsilon\in[0,\frac14]$, the system $(T_\epsilon,I^\Lambda)$ has a unique
  absolutely continuous invariant probability measure
  $d\mu_{\epsilon,\Lambda}=h_{\epsilon,\Lambda}dm^\Lambda$ with the following
  properties:
  \begin{enumerate}[i)]
  \item $(T_\epsilon,\mu_{\epsilon,\Lambda})$ is mixing. Indeed, it has
    exponentially decreasing correlations in time for smooth observables (with
    a speed of decay that depends heavily on the system size $|\Lambda|$,
    though).
  \item The density $h_{\epsilon,\Lambda}$ is of bounded variation, so in
    particular $h_{\epsilon,\Lambda}\in L^{1+1/(|\Lambda|-1)}_{m^\Lambda}$, and
    $h_{\epsilon,\Lambda}>0$ Lebesgue-almost everywhere.
  \end{enumerate}
  The same results remain true if instead of periodic boundary conditions one
  prescribes fixed or free boundary conditions.
\end{thmm}
We give the proof of this theorem in section~\ref{sec:phasetr-bif}. It is a
simple variant of the folklore type argument from \cite{Keller-1997}.

Finally, in section~\ref{sec:smooth}, we discuss how the same results as above
can be produced for an analytic expanding circle map $\tau$. We are not able,
however, to replace also the diffusive nearest neighbour coupling by a
coupling that can be described by a diffeomorphism of the state space
$\Omega$, which on the one hand is close enough to the identity for the
existence part of Theorem~\ref{theo:main} but on the other hand is
sufficiently far away from the identity to admit a phase transition. So we are
not able to construct examples of phase transitions in a class of systems as
studied \eg in
\cite{BEIJK-98,Bardet-02,BK-95,BK-96,BS-88,FR-00,Ji-03,Rugh-02}.

\section{The local map}
\label{sec:local-map}

The basic ingredient of the construction is a continuous piecewise linear
Markov map $\tilde\tau:[-1,1]\to[-v,v]$ where $0<v<1$ and (for later 
use) $0<a<b<c<d<u<v$ are suitable numbers. The map $\tilde\tau$ is symmetric
in the sense that $\tilde\tau(-x)=-\tilde\tau(x)$, it leaves the two intervals
$[-v,-c]$ and $[c,v]$ invariant, and its restrictions to each of these
intervals are mixing with a strictly positive invariant density. However, the
invariant density on $[c,v]$, call it ${\tilde h}$, is highly concentrated on
the subinterval $[d,u]\subset[c,v]$. The map $\tau$ on $[-1,1]$ is then
defined as
\begin{equation}
  \label{eq:tau-def}
  \tau(x)=\hat\tau^3(x)\text{\; where\; }
  \hat\tau(x)=\frac1v\cdot \tilde\tau^k(x)
\end{equation}
with a suitable (rather large) $k\in\NN$.
In particular, $k$ will be
chosen such that $|\hat\tau'|\geq4$.

\begin{figure}
\begin{picture}(500,0)%
\includegraphics{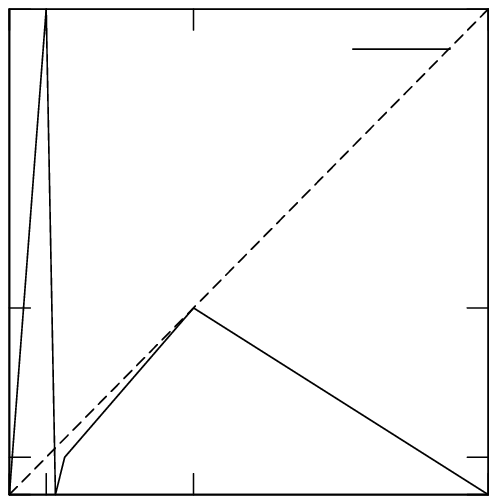}%
\includegraphics{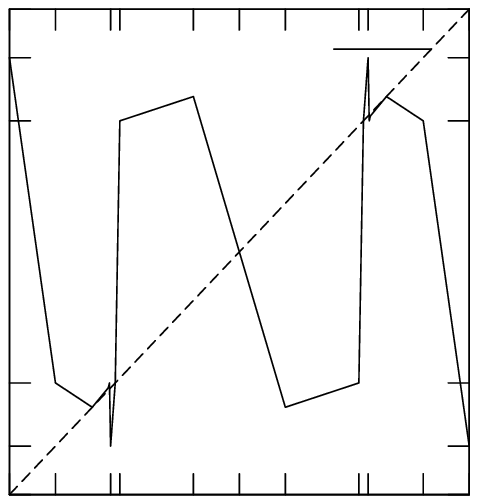}%
\end{picture}%
\begingroup
\setlength{\unitlength}{0.0200bp}%
\begin{picture}(8820,8640)(24900,0)%
\put(825,8090){\makebox(0,0)[r]{\strut{}v}}%
\put(825,3788){\makebox(0,0)[r]{\strut{}u}}%
\put(825,1638){\makebox(0,0)[r]{\strut{}d}}%
\put(825,1100){\makebox(0,0)[r]{\strut{}c}}%
\put(7995,550){\makebox(0,0){\strut{}v}}%
\put(3752,550){\makebox(0,0){\strut{}u}}%
\put(1630,550){\makebox(0,0){\strut{}d}}%
\put(1100,550){\makebox(0,0){\strut{}c}}%
\put(5770,7515){\makebox(0,0)[r]{\strut{}$\tilde\tau$}}%
\end{picture}%
\setlength{\unitlength}{0.0200bp}%
\begin{picture}(8820,8640)(24900,0)%
\put(1100,8090){\makebox(0,0)[r]{\strut{}1}}%
\put(1100,6482){\makebox(0,0)[r]{\strut{}c}}%
\put(1100,7391){\makebox(0,0)[r]{\strut{}v}}%
\put(1100,1799){\makebox(0,0)[r]{\strut{}-v}}%
\put(1100,2708){\makebox(0,0)[r]{\strut{}-c}}%
\put(1100,1100){\makebox(0,0)[r]{\strut{}-1}}%
\put(7995,550){\makebox(0,0){\strut{}1}}%
\put(7333,550){\makebox(0,0){\strut{}v}}%
\put(6639,550){\makebox(0,0){\strut{}d}}%
\put(6306,550){\makebox(0,0){\strut{}b}}%
\put(5347,550){\makebox(0,0){\strut{}a}}%
\put(4685,550){\makebox(0,0){\strut{}0}}%
\put(4023,550){\makebox(0,0){\strut{}-a}}%
\put(3064,550){\makebox(0,0){\strut{}-b}}%
\put(2631,550){\makebox(0,0){\strut{}-d}}%
\put(2037,550){\makebox(0,0){\strut{}-v}}%
\put(1375,550){\makebox(0,0){\strut{}-1}}%
\put(5770,7515){\makebox(0,0)[r]{\strut{}$\tilde\tau$}}%
\end{picture}%
\endgroup
  \caption{The graphs of $\tilde\tau$ (l.h.s.) and of $\tilde\tau|_{[c,d]}$
    (r.h.s.) with parameters $\eta=\frac15$ and $\delta=\frac1{50}$. The
    parameter $\gamma$ is set to $\frac1{200}$ for this figure.  It has to be
    chosen much smaller in order to construct a phase transition.}
  \label{fig:f}
\end{figure}

The construction of the map $\tilde\tau$ depends on two parameters
$\delta,\eta>0$ which obey the inequalities
\begin{equation}
  \label{eq:delta-eta}
  \eta<\frac14\qquad\text{and}\qquad
  \frac{\eta^3}{2 - 4 \eta} < \delta < \frac{\eta^3 - 3 \eta^2 +
    \eta}{4 - 2 \eta}
\end{equation}
\begin{figure}
\begin{picture}(250,125)%
\includegraphics{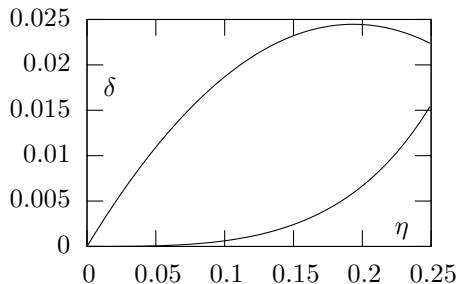}%
\end{picture}%
\begingroup
\setlength{\unitlength}{0.0200bp}%
\begin{picture}(0,0)(12500,0)%
\put(2750,1650){\makebox(0,0)[r]{\strut{} 0}}%
\put(2750,2506){\makebox(0,0)[r]{\strut{} 0.005}}%
\put(2750,3362){\makebox(0,0)[r]{\strut{} 0.01}}%
\put(2750,4218){\makebox(0,0)[r]{\strut{} 0.015}}%
\put(2750,5074){\makebox(0,0)[r]{\strut{} 0.02}}%
\put(2750,5930){\makebox(0,0)[r]{\strut{} 0.025}}%
\put(3025,1100){\makebox(0,0){\strut{} 0}}%
\put(4322,1100){\makebox(0,0){\strut{} 0.05}}%
\put(5619,1100){\makebox(0,0){\strut{} 0.1}}%
\put(6916,1100){\makebox(0,0){\strut{} 0.15}}
\put(8213,1100){\makebox(0,0){\strut{} 0.2}}%
\put(9509,1100){\makebox(0,0){\strut{} 0.25}}%
\put(3500,4600){\makebox(0,0){\strut{}$\delta$}}%
\put(9000,2000){\makebox(0,0){\strut{}$\eta$}}%
\end{picture}%
\endgroup 
\caption{\label{fig:eta-delta}Admissible $(\eta,\delta)$-pairs: the region
  between the two curves.}
\end{figure}
The range of
admissible values of $\eta$ and $\delta$ is the region between the two curves
in Figure~\ref{fig:eta-delta}. One possible choice is $\eta=\frac1{5}$ and
$\delta=\frac1{50}$.  So we use these values in the sequel to illustrate our
construction, but each other choice of $\eta$ and $\delta$ satisfying
\eqref{eq:delta-eta} would do as well.  Figure~\ref{fig:f} shows the graphs of
$\tilde\tau$ and of $\tilde\tau|_{[c,v]}$ that will be defined using these
parameters. Indeed, let
\begin{eqnarray*}
  &&a=1-4\eta\\
  &<&b=1-2\eta-4\delta\\
  &<&c=1-2\eta-3\delta\\
  &<&d=1-2\eta-2\delta\\
  &<&u=1-2\eta+\eta^2=v^2\\
  &<&v=1-\eta
\end{eqnarray*}
and, with a suitable (rather small) $\gamma>0$,
\begin{equation}
  \label{eq:gamma}
  \begin{split}
  d\quad<&\quad d'=d+\gamma\\
  <&\quad d''=d+2\gamma\\
  <&\quad u\ .    
  \end{split}
\end{equation}
Observe that, by this choice, 
\begin{equation}
  \label{eq:c2<v3}
  c<u=v^2\;.
\end{equation}
Define $\tilde\tau|_{[0,1]}$ as the piecewise linear interpolation of
\begin{quote}
  $\tilde\tau(0)=0$, $\tilde\tau(a)=-u$, $\tilde\tau(b)=-c$,
  $\tilde\tau(c)=c$, $\tilde\tau(d)=v$, $\tilde\tau(d')=c$,
  $\tilde\tau(d'')=d$, $\tilde\tau(u)=u$, $\tilde\tau(v)=c$,
  $\tilde\tau(1)=-v$,
\end{quote}
and extend $\tilde\tau$ to all of $[-1,1]$ by
$\tilde\tau(-x)=-\tilde\tau(x)$. 

By definition, $\tilde\tau$ and also the restriction $\tilde\tau|_{[c,v]}$ are
piecewise linear Markov maps, so their invariant densities can be calculated
as eigenvectors of corresponding Markov matrices.
The unique $\tilde\tau$-invariant probability
density ${\tilde h}$ on $[c,v]$ is constant
on $[c,d),[d,u]$ and on $(u,v]$.  Using a computer algebra system one checks
the following expressions for the corresponding probabilities under ${\tilde
  h}(x)dx$:
\begin{equation}
  \label{eq:h0prob}
  \begin{split}
    \int_c^d{\tilde h}(x)\,dx &= 
    \frac{\left( \eta^2 + 4 \delta \right)
    }{\eta^4 + 4 \delta \eta^2 + 4 \delta^2}\,\gamma + {\rm O}(\gamma^2) 
    \\
    \int_d^u{\tilde h}(x)\,dx &= 1 -\frac{3\eta^3+10\delta\eta+12\delta^2-3\delta\eta^2-2\eta^4}{\eta^5 + 3 \delta \eta^4 + 4
      \delta \eta^3 + 12 \delta^2 \eta^2 + 4 \delta^2 \eta + 12
      \delta^3}\,\gamma + {\rm O}(\gamma^2)
    \\
    \int_u^v{\tilde h}(x)\,dx &=  \frac{6\delta\eta+2\eta^3-6\delta\eta^2-2\eta^4}{\eta^5 + 3 \delta \eta^4 + 4
      \delta \eta^3 + 12 \delta^2 \eta^2 + 4 \delta^2 \eta + 12
      \delta^3}\,\gamma + {\rm O}(\gamma^2)
      \;.
  \end{split}
\end{equation}
(For $\eta=\frac15$ and $\delta=\frac1{50}$ the corresponding values are
approximately: $18.75\gamma+{\rm O}(\gamma^2)$, $1 - 37.98\gamma+ {\rm
  O}(\gamma^2)$ and $19.23\gamma+{\rm O}(\gamma^2)$.)  Furthermore, it is easy
to check that each nontrivial interval contains a point which, under iteration
of $\tilde\tau$, eventually is mapped into $[-v,-c]\cup[c,v]$. Hence, all
invariant probability
densities of $\tilde\tau$ are convex combinations 
\begin{equation}
  \label{eq:tilde-h-convex}
  \tilde h_\alpha(x):=\alpha\,\tilde h(x)+(1-\alpha)\,\tilde h(-x)
  \quad(0\leq\alpha\leq1)\;.
\end{equation}
Because of \eqref{eq:h0prob}, for each admissible pair $(\eta,\delta)$ and for
each $\kappa>0$ one can choose $\gamma>0$ sufficiently small such that, for
all $\alpha\in[0,1]$,
\begin{equation}
  \label{eq:halpha-prob}
  \int_{[-u,-d]\cup[d,u]}\tilde h_\alpha\,dm>1-\kappa\;.
\end{equation}

Let $m$ be Lebesgue measure and denote by $P_\tau,
P_{\tilde\tau}$ and $P_{\hat\tau}$ the Perron-Frobenius operators of $\tau,
\tilde\tau$ and $\hat\tau$, respectively. As $\tilde\tau|_{[c,v]}$ is mixing
and because of \eqref{eq:tilde-h-convex}, there are $C>0$ and $\rho\in(0,1)$
with the following property:
for each $f\in L^1_m(I)$ there is $\alpha\in[0,1]$ such that, for all $k>0$,
\begin{equation}
  \label{eq:tilde-tau-mixing}
  \int\left|(P_{\tilde\tau}^kf)-\int f\,dm\cdot{\tilde h}_\alpha\right|dm
  \leq
  C\,\rho^k\,\Var(f)\;.
\end{equation}
 see \eg \cite[Theorem 2.1]{KeLi-2005a}.

To proceed we need the following family of probability densities on $I$:
\begin{equation}
  \label{eq:hat-h-alpha}
  \hat h_\alpha(x):=
  v\cdot\tilde h_\alpha(vx)\quad(0\leq\alpha\leq1)\;.
\end{equation}
The following two lemmas collect all the information on $\hat\tau$ and
$P_{\hat\tau}$ we need.

\begin{lem}
  \label{lem:single-site-1}
  The map $\hat\tau=\frac 1 v \tilde\tau^k:[-1,1]\to[-1,1]$ is surjective and
  mixing. Indeed, if $J$ is any maximal monotonicity interval of $\hat\tau$,
  then $\hat\tau^3(J)=[-1,1]$.
\end{lem}
\begin{proof}
  $\hat\tau$ is surjective because $\tilde\tau([-v,v])=[-v,v]$. In order
  to prove that the piecewise expanding map $\hat\tau$
 is mixing, it suffices to show that
  $\hat\tau^3(J)=[-1,1]$ for each maximal monotonicity interval $J$ of
  $\hat\tau$, see \cite{Bowen-1977}: If $J$ is a maximal monotonicity interval
  of $\hat\tau$, then it is also a maximal monotonicity interval of
  $\tilde\tau^k$. As
  \begin{displaymath}
    -1<-v<-u<-d'<-d<-c<-a<a<c<d<d'<u<v<1
  \end{displaymath}
  defines a Markov partition for $\tilde\tau$ (all branches with respect to
  this partition are monotone, but
  not all are linear neither are all these branches maximal monotone!),
  $\tilde\tau^kJ$ contains
one of the images of the Markov intervals, \ie one of the
  intervals
  \begin{displaymath}
    [-c,v],[-u,-c],[-v,-c],[-c,u],[-u,u],[-u,c],[c,v],[c,u]\;.
  \end{displaymath} 
  Each of these intervals contains at least one of $[-u,-c]$ or $[c,u]$, so we
  will assume from now on that $\tilde\tau^k J\supseteq[c,u]$, the other case
  being treated in the same way because of the symmetry of $\tilde\tau$.  Then
  $\hat\tau J\supseteq[\frac cv,\frac uv]$. As $\frac uv=v$ and
  $\tilde\tau(v)=c$, we see that $\tilde\tau\hat\tau J$ contains an interval
  $[c,c+q]$ of length $q=\frac{u-c}v\cdot\frac{u-c}{v-u}
  =\frac{(\eta^2+3\delta)^2}{\eta(1-\eta)^2}$. Choose $k$ big enough such that
  this interval is mapped by $\tilde\tau^{k-2}$ at least over the interval
  $[c,d]$. (For $\eta=\frac15$ and $\delta=\frac1{50}$ already $k\geq2$ would
  do.) Then $\tilde\tau^k\hat\tau J\supseteq[c,v]$, so $\hat\tau^2
  J\supseteq[\frac cv,1]\supset[v,1]$. Hence
  $\tilde\tau\hat\tau^2J\supseteq[-v,c]$, so
  $\tilde\tau^k\hat\tau^2\supseteq[-v,v]$ and thus
  $\hat\tau^3\supseteq[-1,1]$.
\end{proof}
\begin{lem}
  \label{lem:single-site-2}Recall that $\hat\tau=\frac1v\tilde\tau^k$ and
  denote by $\Ts:I^s\to I^s$ the $s$-fold direct product of $\hat\tau$ with
  itself. (We will be using the case $s=4$ later.)
  \begin{enumerate}[a)]
  \item\label{lem:single-site-2-a} There is a constant ${\beta}>0$ (whose choice
    depends only on $\eta,\delta$ and $\gamma$) such that, for each
    $\kappa\in(0,\frac12)$ and sufficiently large $k$, the following estimate
    holds for each $f\in L^1_{m}(I^s)$:
    \begin{equation}
      \label{eq:LY-1D}
      \var(P_{\Ts}f)
      \leq
      \kappa\,\var(f)+{\beta}\,\int|f|\,dm\,.
    \end{equation}
  \item\label{lem:single-site-2-b} Given $\kappa>0$, the following holds for
    all sufficiently large $k$: for each $f\in L^1_{m}(I^s)$ of bounded
    variation there exists $\alpha=(\alpha_1,\dots,\alpha_s)\in[0,1]^s$ such
    that
    \begin{equation}
      \label{eq:JB1}
      \var\left(P_{\Ts}f-\int f\,dm\cdot\hat h_\alpha\right)
      \leq
      \kappa\,\var(f)
    \end{equation}
    where $\hat h_\alpha(x_1,\dots,x_s):=\hat
    h_{\alpha_1}(x_1)\cdot\dots\cdot\hat h_{\alpha_s}(x_s)$.
  \item\label{lem:single-site-2-c} Given $\kappa>0$, for all sufficiently
    small $\gamma$ and all $\alpha\in[0,1]^s$,
    \begin{equation}
      \label{eq:kappa}
      \int_{([-\frac uv,-\frac dv]\cup[\frac dv,\frac uv])^s}{\hat
        h_\alpha}\,dm
      >
      1-\kappa\;.
    \end{equation}
  \end{enumerate}
\end{lem}
\begin{proof}
  \ref{lem:single-site-2-a})\;\eqref{eq:LY-1D} is the finite-dimensional
  (uncoupled) Lasota-Yorke inequality, see \eg \cite[Lemma 3.2 and eq.
  (60)]{KeLi-2005a}.  The other two estimates are immediate consequences.
  \\[2mm]
  \ref{lem:single-site-2-b})\; This follows easily from
  \eqref{eq:tilde-tau-mixing} and \eqref{eq:hat-h-alpha}, see \eg
  \cite[section 4.5]{KeLi-2005a}.
  \\[2mm]
  \ref{lem:single-site-2-c})\; This is an immediate consequence of
  \eqref{eq:hat-h-alpha} and \eqref{eq:halpha-prob} (where the $\kappa$ from
  \eqref{eq:kappa} equals $1-(1-\kappa)^s$ in terms of the $\kappa$ from
  \eqref{eq:halpha-prob}).
\end{proof}

Hence, the ``local units'' $\tau=\hat\tau^3$ will behave like Markov chains
that admit transitions from the positive to the negative part of $[-1,1]$ and
vice versa only with very small probabilities because, for sufficiently large
$k$, most of the mass will be concentrated on $[-\frac uv,-\frac dv]\cup[\frac
dv,\frac uv]\subset[-v,-c]\cup[c,v]$ after each application of $\hat\tau$, and
points in this set will not change sign during the next application of
$\hat\tau$.  Nevertheless, there is always a small mass, proportional to
$\gamma$, which is mapped under the next application of $\tilde\tau$ to
$[-v,-c]$, so transitions are never completely excluded.

In the next section we will see how coupling with at least one neighbour of
the same sign can prevent such a sign flip, and we will be careful enough to
make sure that coupling with two neighbours of opposite sign indeed forces a
flip.

\begin{rem}\label{rem:full-branches}
  In the rest of this section and in
  sections~\ref{sec:local-coupling}~-~\ref{sec:exp-estimate}, only those
  properties of $\hat\tau$ will be used that are formulated in
  Lemma~\ref{lem:single-site-2}. (Lemma~\ref{lem:single-site-1} will play a
  role only in the proof of Theorem~\ref{theo:finite-dim} in
  section~\ref{sec:phasetr-bif}.) Therefore, the same proof of
  Theorem~\ref{theo:main} will work for two kinds of modifications of
  $\hat\tau$.
  \begin{enumerate}[(A)]
  \item\label{item:increasing} All arguments remain unchanged, if we replace
    the decreasing branches of $\tilde\tau$ by increasing ones with the same
    domain and the same range: we obtain a piecewise linear Markov map
    $\acute\tau$, it satisfies the Lasota-Yorke inequality with the same
    constants, and it has the same associated Markov chain as $\tilde\tau$. So
    it has in particular the same 2-dimensional space of invariant densities
    $\tilde h_\alpha$. $\acute\tau$ is of course no longer continuous, but the
    continuity of $\tilde\tau$ was not used in the proof of 
    Lemma~\ref{lem:single-site-2}.
  \item\label{item:smooth} Let $\check\tau=\frac1v\acute\tau^k$. We may also
    replace $\check\tau$ by any map $\bar\tau$ which has full branches and
    whose Perron-Frobenius operator is a small perturbation of that of
    $\check\tau$ in the sense that it satisfies a Lasota-Yorke inequality with
    the same constants as $\hat\tau$ and $\check\tau$ (thus leading to the
    same choices of ${\beta}$ and $k=k(\kappa)$ in
    Lemma~\ref{lem:single-site-2}\ref{lem:single-site-2-a}) and that
    \begin{displaymath}
      |||P_{\bar\tau}-P_{\check\tau}|||
      :=
      \sup\left\{\int|(P_{\bar\tau}-P_{\check\tau})(f)|\,dm:\ 
        \Var(f)\leq1\right\}
    \end{displaymath}
    is sufficiently small. For such maps, the conclusions of
    Lemma~\ref{lem:single-site-2}\ref{lem:single-site-2-a} hold trivially, and
    the persistence of the conclusions of
    Lemma~\ref{lem:single-site-2}\ref{lem:single-site-2-b}/\ref{lem:single-site-2-c}
    is a consequence of the spectral stability theorem from \cite{KeLi-1998},
    see also \cite[Sections 3.1 - 3.3]{Baladi-book} for a coherent account of
    the spectral (perturbation) theory for Perron-Frobenius operators of
    piecewise expanding transformations.
  \end{enumerate}
  In section~\ref{sec:smooth} we will
  come back to these remarks.
\end{rem}

\section{The coupling: local effects}
\label{sec:local-coupling}
The coupling map $\Phi_\epsilon:\Omega\to\Omega$ can be described in terms of a
local coupling rule $\phi_\epsilon:I^3\to I$,
\begin{equation}
  \begin{split}
  \label{eq:phi-def}
  \phi_\epsilon(x,y,z):=&(1-\epsilon) x+\frac\epsilon2(y+z)\;,\text{\; so that  }\\
  (\Phi_\epsilon(\xx))_\ii=&\phi_\epsilon(\xx_\ii,\xx_{\ii+\ee_1},\xx_{\ii+\ee_2})\;.
  \end{split}
\end{equation}
Suppose now that $x,y,z\in[-1,1]$, $\frac dv\leq|x|,|y|,|z|\leq\frac uv=v$,
and $\epsilon\in[0,\eta]$. Denote $x'=\phi_\epsilon(x,y,z)$.  By symmetry we
may assume without loss of generality that $x>0$.  Our first observation is
obvious,
  \begin{equation}
    \label{eq:case1}
    x'
    \leq
    v\;.
  \end{equation}
  \vspace*{0mm}
\begin{description}
\item [\hspace*{4mm}Suppose that $y>0$ or $z>0$] Without loss of
  generality let $y>0$. Then
  \begin{equation}
    \label{eq:case2}
    \begin{split}
      x'
      &\geq
      (1-\eta)\frac dv+\frac\eta2\,\frac{d-u}v
      =
      d+\frac{\eta(\eta^2+2\delta)}{2(1-\eta)}
      > c
    \end{split}
  \end{equation}
  where we used $\delta>\frac{\eta^3}{2-4\eta}$ for the last inequality,
  see \eqref{eq:delta-eta}.
  \vspace*{2mm}
\item [\hspace*{4mm}Suppose that $y,z<0$ and $\epsilon=\eta$] Then
  \begin{equation}
    \label{eq:case3a}
    x'
    \geq
    (1-\eta)\frac dv-\eta v
    =
    d-\eta v
    >
    a
  \end{equation}
  where we used for the last equation that $\delta < \frac{\eta^3 - 3 \eta^2 +
    \eta}{4 - 2 \eta}\leq\frac{\eta+\eta^2}2$ (for $\eta\in[0,1]$), see also
  \eqref{eq:delta-eta}. In this case we also have
  \begin{equation}
    \label{eq:case3b}
    \begin{split}
      x'
      &\leq
      (1-\eta)v-\eta\frac dv
      =
      (1-\eta)^2-\frac{\eta(1-2\eta-2\delta)}{1-\eta}
      <
      b
    \end{split}
  \end{equation}
  where we used again $\delta < \frac{\eta^3 - 3 \eta^2 +
    \eta}{4 - 2 \eta}$ for the last
  inequality.
  We finally notice that inequalities \eqref{eq:case3a} and \eqref{eq:case3b} are strict, hence still satisfied for values of $\epsilon$ which are smaller than $\eta$, but close to it, say for $\epsilon\in[\epsilon_2,\eta]$.
  \vspace*{2mm}
\end{description}
We summarize these observations: 
\begin{itemize}
\item If at least one of $y$ and $z$ is in $[\frac dv,v]$, then $x'\in(c,v]$
  and hence $\tilde\tau^k x'\in[c,v]$.
\item If both $y$ and $z$ are in $[-v,-\frac dv]$ and if $\epsilon\in[\epsilon_2,\eta]$, then
  $x'\in(a,b)$ so that $\tilde\tau(x')\in[-u,-c]$ and hence
  $\tilde\tau^kx'\in[-v,-c]$.

\end{itemize}

Indeed, we will see that, with very high
probability, $\frac dv\leq|\hat\tau^3(x')|\leq v$ again so that, if $x$
is the state of the coupled system at a given time at lattice site $\ii$ and
if $y$ and $z$ are the states at the same times at sites $\ii+\ee_1$ and
$\ii+\ee_2$, then the sign of $x'$, which is the updated value at site $\ii$,
deviates only with small probability from the ``majority vote'' of $x,y$ and
$z$, i.e. from $\sign(\sign x+\sign y+\sign z)$. Hence this interaction mimics
the local rule of Toom's probabilistic cellular automaton \cite{Toom-1980}. In
the next section we will see how the proof that this automaton admits at least
two invariant measures can be modified for our purpose. We follow the
treatment of \cite{LeMaSp-1990}.

\section{The coupling: global effects}
\label{sec:global-coupling}

It is useful to introduce the usual total variation norm on signed measures:
\begin{equation}
  \label{eq:regularnorm}
  |\mu|:=\sup_{|\vf|_{\Co^0(\Omega)}\leq 1}\mu(\vf)\;.
\end{equation}
Just like in \cite[Sect.~3.3]{KeLi-2005a} one checks easily that
\begin{equation}
  \label{eq:two-norms}
|\mu|\leq\frac12\Var\mu\;.
\end{equation}

We consider the
dynamics acting directly on the measures via the linear operator
$T^{\ast}_\eps\mu(A):=\mu(T_\eps^{-1}A)$ (for each measurable set $A$).
The basic facts concerning the operator $T^{\ast}_\eps$ are detailed in the
following lemma.
\begin{lem}[Lasota-Yorke inequality]\label{lem:lasota-yorke}
  Recall that $\hat\tau=(\frac1v \tilde\tau^k)$.  There exists constants
  $k_0>0$ and $B>0$ which depend only on the parameters $\eta,\delta$ and
  $\gamma$ such that, for all integers $k\geq k_0$, $|\hat\tau'|\geq4$ and
  such that, for each $|\eps|\leq\frac14$, the operator $T^{\ast}_\eps$ is
  well defined as an operator on $\BV$. In addition, for each $\mu\in\BV$
  holds true
  \[
  \begin{split}
    |T^{\ast}_\eps\mu|&\leq |\mu|\\
    \Var (T^{\ast n}_\eps\mu)&\leq 2^{-n}\Var\mu+B|\mu|\;.
  \end{split}
  \]
  In particular, $\Var(T^\ast_\eps\mu)\leq 2B$ for each probability measure
  $\mu\in\BV$ with $\Var\mu\leq2B$, and $\Var\mu\leq B$ for each
  $T_\epsilon$-invariant probability measure $\mu\in\BV$.
  
  The same assertions hold for the case of a finite lattice
  $\Lambda=(\ZZ/L\ZZ)^2$.
\end{lem}
\begin{proof}
  This follows from the special case $\theta=1$ of Proposition 4.1 in
  \cite{KeLi-2005a}. (For finite $\Lambda$ use \cite[Proposition
  3.2]{KeLi-2005a}.)  Observe that the proof given there for $\Lambda=\ZZ$
  applies (only if $\theta=1$!) without changes to $\Lambda=\ZZ^2$.  The
  particular constants in the present lemma follow from the proofs in
  \cite[section 3]{KeLi-2005a} by observing the following facts:
  \begin{enumerate}[$\triangleright$]
  \item Our coupling is a $(1,0)$-coupling in the sense of
    \cite[section 3.1]{KeLi-2005a}.
  \item We consider the local map $\tau$ itself, so $\ell=1$ in
    \cite[section 3.4]{KeLi-2005a}.
  \item Recall that $\tau(x)=\hat\tau^3(x)=(\frac1v
    \tilde\tau^k)^3(x)$. We choose $k$ so large that
    $|\hat\tau'|\geq4$.  An elementary detailed analysis of the map
    shows that this is always possible because
    $|(\tilde\tau^2|_{[c,v]})'|\geq\min\{9,(\frac{2\delta+\eta^2}{2\delta+\eta^2-3\gamma})^2,\frac{3\delta}\gamma\}$.
    Indeed, since $\gamma$ is rather small, the most critical slope of
    $\hat\tau$ occurs in a small one-sided left neighbourhood of $-a$,
    which is mapped to a one-sided left neighbourhood of the fixed
    point $u$. This slope is
    $\frac{\eta^2+3\delta}{2\eta-4\delta}(\frac{2\delta+\eta^2}{2\delta+\eta^2-3\gamma})^{k-1}$.
    With these choices, Lemmas~3.2 and 3.3 in \cite{KeLi-2005a} yield the
    claimed values for the constants when choosing $\epsilon_1=\frac14$ there.
    Observe also that the exponential factor $2^{-n}$ can be replaced by any
    factor $\rho^n$, $\rho\in(0,1)$, at the price of larger $k$ and $B$.
  \end{enumerate}
\end{proof}
\begin{rem} 
  The choice of $k$ we make in Lemma \ref{lem:lasota-yorke} is sufficient for
  the proof of Theorem~\ref{theo:main}. For Theorem~\ref{theo:finite-dim}, our
  proof will require to choose a still larger $k$~: we require that
  $|\hat\tau'|\geq\frac{12v}{u-c}$, which implies in
  particular that $|\hat\tau'|\geq4$.
\end{rem}
\noindent
\textbf{Proof of Theorem~\ref{theo:main}\ref{theo:main-a}}\; The existence of
invariant probability measures in $\BV$ follows from
Lemma~\ref{lem:lasota-yorke} just as in \cite[section 4.4]{KeLi-2005a}. The
proof given there for the lattice $\Lambda=\ZZ$ applies without changes to the
lattice $\Lambda=\ZZ^2$.
\quad\\[3mm]
\noindent
\textbf{Proof of Theorem~\ref{theo:main}\ref{theo:main-b}}\;
This follows at once from \cite{KeLi-2005b}.
\\[3mm]
\noindent
\textbf{Proof of Theorem~\ref{theo:main}\ref{theo:main-c}}\; 
In order to get close to the formal setting of \cite{LeMaSp-1990}, we
introduce the short hand notation
\begin{equation}
  \label{eq:short-hand}
  \xx^t:=T_\epsilon^t(\xx)\hspace*{5mm}\text{for $\xx\in\Omega$ and $t\in\NN$.}
\end{equation}
In \cite{LeMaSp-1990}, the authors work with time indexed by $t\in\ZZ$, in
particular they specify certain events at times $-N$ and $0$. In our setting
it seems more natural to shift these events to times $0$ and $N$, and since
all estimates are for fixed $N$, this is just a matter of convenience.

In order to make clear which estimates we need, we reproduce here a sketch of
the basic argument in \cite{LeMaSp-1990} with notations adapted to our
setting.  Fix $N\in\NN$ and let
\begin{displaymath}
  \Sigma_+^N
  :=
  \{{\bsigma}=(\sigma_\ii^t)\in\{-1,+1\}^{\Lambda\times\{0,\dots,N\}}:
  \sigma_\nn^N=-1\text{ and }\sigma_\ii^0=+1\,\forall\ii\in\Lambda\}
\end{displaymath}
In \cite[Appendix A]{LeMaSp-1990} the authors construct a family $\V^{(N)}$ of
nonempty finite subsets $\hat V\subset\Lambda\times\NN$ together with a map
$V^{(N)}:\Sigma_+^N\to\V^{(N)}$ in such a way that\,\footnote{More correctly,
  what we state here is a consequence of what is proved in \cite{LeMaSp-1990}.
  It is precisely what we need for our proof, and a reader who attempts to go
  through the proof in \cite{LeMaSp-1990} will find it a good warm-up exercise
  to check that our claims (i) and (ii) are indeed immediately implied by what
  is proved there.}
\begin{enumerate}[(i)]
\item $\card{\V_M^{(N)}}\leq (48)^{8M}$ for each $M\geq1$ where
  $\V_M^{(N)}:=\{\hat V\in\V^{(N)}:\card\hat V=M\}$.
\item Let $\bsigma\in\Sigma_+^N$. Each $(\ii,t)\in V^{(N)}(\bsigma)$ is an
  \emph{error site} in space-time for the configuration $\bsigma$, \ie
  $\sigma_\ii^{t+1}\neq\sign(\sigma_\ii^t+\sigma_{\ii+\ee_1}^t+\sigma_{\ii+\ee_2}^t)$.
  In other words, $\bsigma$ deviates from Toom's deterministic majority rule at
  $(\ii,t)$.
\end{enumerate}

Now we apply this purely combinatorial construction to our problem.
Let 
\begin{displaymath}
  \Omega_+
  :=
  \{\xx\in\Omega:\,x_\ii^0=x_\ii>0\,\forall\ii\in\Lambda\}\text{ and }
  \Omega_+^N
  :=
  \{\xx\in\Omega_+:\,x_\nn^N\leq0\}\;.
\end{displaymath}
To each $\xx\in\Omega_+^N$ we associate $\bsigma(\xx)\in\Sigma_+^N$ defined by
$\bsigma_\ii^t=+1$ or $-1$ according to whether $x_\ii^t>0$ or $x_\ii^t\leq0$,
respectively. (We do not claim that each $\bsigma\in\Sigma_+^N$ occurs as some
$\bsigma(\xx)$.) For a finite subset $W$ of $\Lambda\times\NN$ let
\begin{equation}
  \label{eq:error-sites}
  E(W):=\{\xx\in\Omega_+:\,\text{Each $(\ii,t)\in W$
      is an error site for $\bsigma(\xx)$}\}\;.
\end{equation}
Then, if $\mu$ is any probability measure on $\Omega_+$,
\begin{equation}
  \label{eq:LMS-basic}
  \begin{split}
    \mu(\Omega_+^N)
    &=
    \sum_{M=1}^\infty\sum_{\hat V\in\V_M^{(N)}}
    \mu\{\xx\in\Omega_+^N:\,V^{(N)}(\bsigma(\xx))=\hat V\}\\
    &\leq
    \sum_{M=1}^\infty\sum_{\hat V\in\V_M^{(N)}}
    \mu(E(\hat V))\\
    &\leq
    \sum_{M=1}^\infty (48)^{8M}\cdot q_M
  \end{split}
\end{equation}
where 
\begin{equation}
  \label{eq:q_M}
  q_M:=\sup\{\mu(E(\hat V)): N\geq0\text{ and }\hat V\in\V_M^{(N)}\}\;.
\end{equation}

In the next section we will derive the exponential estimate
\begin{equation}
  \label{eq:exp-estimate}
  \text{$q_M\leq\Delta^{M}$ with $\Delta=\frac14(48)^{-8}>0$}
\end{equation}
valid when $\mu=\lambda_+$ is the product Lebesgue measure
on $[\frac dv,v]^\Lambda$.
Hence, for all $N>0$,
\begin{equation}
  \label{eq:final-estimate}
  (T_\epsilon^{*N}\lambda_+)\{\xx\in\Omega:x_\nn\leq0\}
  =
  \lambda_+(\Omega_+^N)\leq\frac13\;.
\end{equation}
Now, as in \cite[proof of Theorem 4.1]{KeLi-2005a} it follows that
$\frac1n\sum_{N=1}^nT_\epsilon^{*N}\lambda_+$ converges weakly to a
$T_\epsilon$-invariant probability measure $\mu_+\in\BV$. So
$\mu_+\{\xx\in\Omega:x_\nn\leq0\}\leq\frac13$. Interchanging the roles of
$+$ and $-$ one also finds a $T_\epsilon$-invariant probability $\mu_-\in\BV$
with $\mu_-\{\xx\in\Omega:x_\nn\geq0\}\leq\frac13$.  This finishes the proof
of Theorem~\ref{theo:main}\ref{theo:main-c}, except for the estimate
\eqref{eq:exp-estimate} which is derived in the next section.

\section{The exponential estimate \eqref{eq:exp-estimate}}
\label{sec:exp-estimate}

Let $U(\ii)=\{\ii,\ii+\ee_1,\ii+\ee_2,\ii+\ee_1+\ee_2\}$. Then $x_\ii^{t+1}$
does not depend on $x_\jj^t$ if $\jj\not\in U(\ii)$, and so no $x_\kk^{s}$
influences any $x_\jj^{t}$ as long as $\kk\in U(\ii)$ and
$\jj\in\Pi(\ii)$ where
\begin{equation}
  \label{eq:Pi(i)}
  \Pi(\ii):=\{\jj\in\Lambda:j_1\geq i_1+2\text{ or }j_2\geq i_2+2\}\;.
\end{equation}
Hence $T_\epsilon$ can be ``restricted'' to $I^{\Pi(\ii)}$ as an autonomous
dynamical system, call it $T_{\epsilon,\Pi(\ii)}$, and also to
$I^{U(\ii)\cup \Pi(\ii)}$, call this $T_{\epsilon,U(\ii)\cup \Pi(\ii)}$.
Observe that $T_{\epsilon,U(\ii)\cup \Pi(\ii)}$ is indeed a skew product
transformation over the base $T_{\epsilon,\Pi(\ii)}$. Therefore, the action of
$T_\epsilon$ on the coordinates $x_\kk$, $\kk\in U(\ii)$, can be interpreted
as a nonautonomous dynamical system on $I^{U(\ii)}$ that is governed by the
$T_{\epsilon,\Pi(\ii)}$-orbit of
$\xx_{\Pi(\ii)}$. We use the following notation:
\begin{equation}
  \label{eq:nonautonomous}
  T_{\epsilon,U(\ii)\cup
  \Pi(\ii)}^n(y,\xx_{\Pi(\ii)})=:\Big(T_{\epsilon|\xx_{\Pi(\ii)}}^n(y),T_{\epsilon,\Pi(\ii)}^n(\xx_{\Pi(\ii)})\Big)\;\text{ for }y\in I^{U(\ii)}\text{ and
  }\xx_{\Pi(\ii)}\in I^{\Pi(\ii)}\,,
\end{equation}
and we write $P_{\epsilon|\xx_{\Pi(\ii)}}^n$ for the Perron-Frobenius operator
of $T_{\epsilon|\xx_{\Pi(\ii)}}^n$ on $L^1_{{\rm Leb}}(I^{U(\ii)})$.  This
is the setting studied in \cite{KeZw-2002}. In particular, for
$f\in L^1_{{\rm Leb}}(I^{U(\ii)})$ and each $\xx_{\Pi(\ii)}$,
\begin{equation}
  \label{eq:LY-conditional}
  \var(P_{\epsilon|\xx_{\Pi(\ii)}}^nf)\leq 2^{-n}\var(f)+B\|f\|_1\;,
\end{equation}
where, as in the proof of Lemma~\ref{lem:lasota-yorke}, the factor $2^{-n}$ is
achieved by choosing $k$ in the definition of $\hat\tau=\frac1v\tilde\tau^k$
sufficiently large. 

Now fix some $\hat V\in\V_M^{(N)}$. As
\begin{displaymath}
  \Lambda\times\NN
  =
  (2\Lambda\cup(2\Lambda+\ee_1)\cup(2\Lambda+\ee_2)\cup(2\Lambda+\ee_1+\ee_2))\times(2\NN\cup(2\NN+1))
\end{displaymath}
is the disjoint union of eight sublattices, 
the intersection of $\hat V$ with at least one of these sublattices, call this
intersection $\tilde V$,
satisfies
\begin{displaymath}
  \card{\tilde V}
  \geq
  \frac18\card{\hat V}=\frac M8\;.
\end{displaymath}
We are going to prove that 
\begin{equation}
  \label{eq:exp-est-modif}
  \mu(E(\tilde V))\leq\Delta^{8\cdot\card{\tilde V}}\;.
\end{equation}
As $E(\hat V)\subseteq E(\tilde V)$, estimate \eqref{eq:exp-estimate} then
follows at once.

In order to prove \eqref{eq:exp-est-modif} by induction we introduce, for each
 $\Gamma\subseteq\Lambda$, the set
 \begin{equation}
   \label{eq:VtildeGamma}
   \tilde V_\Gamma:=\tilde V\cap(\Gamma\times\NN)\;.
 \end{equation}
 Denote also, for the moment, $\Lambda':=\{\ii\in\Lambda:\tilde
 V_{\{\ii\}}\neq\emptyset\}$.  By construction of $\tilde V$ from $\hat V$, we
 have $\Lambda'\cap U(\ii)=\Lambda'\cap\{\ii\}$ for each $\ii\in\Lambda'$.  It
 follows that there is (at least) one ``maximal'' site $\ii\in\Lambda'$ in the
 sense that $\Lambda'\setminus\{\ii\}\subseteq\Pi(\ii)$.  We fix such a site
 $\ii$ now. In order to prove
 \eqref{eq:exp-est-modif} it suffices (by induction) to show that
\begin{equation}
  \label{eq:exp-est-induct}
  \mu(E(\tilde V))
  \leq
  \Delta^{8\cdot\card{\tilde V_{\{\ii\}}}}\cdot\mu(E(\tilde V_{\Pi(\ii)}))\;.
\end{equation}


We may write $E(\tilde V)=E(\tilde V_{\{\ii\}})\cap E(\tilde V_{\Pi(\ii)})$,
where $E(\tilde V_{\Pi(\ii)})$ only depends on coordinates at sites from
$\Pi(\ii)$ whereas $E(\tilde V_{\{\ii\}})$ depends on those from
$U(\ii)\cup\Pi(\ii)$.  Hence, for any $\mu\in\ACC$, denoting
$h_{\xx_{\Pi(\ii)}}$ its conditional density on $I^{U(\ii)}$ given
$\xx_{\Pi(\ii)}$ outside, one has
\begin{equation}
\label{eq:space-decoupl}
\begin{split}
  \mu(E(\tilde V)) 
  & = 
  \int_{E(\tilde V_{\Pi(\ii)})}\int 1_{E(\tilde V_{\{\ii\}})}
  (x_{U(\ii)},\xx_{\Pi(\ii)})h_{\xx_{\Pi(\ii)}}(x_{U(\ii)})\,dx_{U(\ii)}\,
  d\mu(\xx_{\Pi(\ii)})\\
  & \leq 
  \sup_{\xx_{\Pi(\ii)}\in I^{\Pi(\ii)}} 
  \Big(\int_{B_0(\xx_{\Pi(\ii)})}h_{\xx_{\Pi(\ii)}}(y)dy\Big)
  \cdot\mu(E(\tilde V_{\Pi(\ii)}))
\end{split}
\end{equation}
with 
\begin{equation}
  \label{eq:B_0}
  B_0(\xx_{\Pi(\ii)})
  =
  \{y\in I^{U(\ii)}\,:\, (y,\xx_{\Pi(\ii)})\in
  E(\tilde V_{\{\ii\}})\}\;.  
\end{equation}


We can now fix $\xx_{\Pi(\ii)}$ and work with the nonautonomous system
$T^t_{\epsilon|\xx_{\Pi(\ii)}}$ on $I^{U(\ii)}$. We also denote, for
$y\in I^{U(\ii)}$,
\begin{equation}
  \label{eq:notationJB}
  \Phi_\epsilon(y,\xx_{\Pi(\ii)}^t)
  =:
  \left(\Phi_{\epsilon}^{\{t\}}(y),
  \Phi_{\epsilon,\Pi(\ii)}(\xx_{\Pi(\ii)}^t)\right)\;,
\end{equation}
see also \eqref{eq:short-hand}, and, with a slight abuse of notation, we denote
by $\hat T$ (resp. $\tilde T$) the 4-fold direct product of $\tilde \tau$
(resp. $\hat\tau$) on $I^{U(\ii)}$.

Let $\xi\in I^{U(\ii)}$. We denote $\xi^t=T^t_{\epsilon|\xx_{\Pi(\ii)}}(\xi)$
and
\begin{equation}
  \label{eq:Mxt}
  \begin{split}
      M(\xi^t):=&\{y\in I^{U(\ii)}:
      \text{
        $y_\ii$ obeys the majority rule relative to $\xi^t$}\}\\
      =&
      \{y\in I^{U(\ii)}:
        \sign(y_\ii)=\sign(\sign(\xi^t_\ii)+\sign(\xi^t_{\ii+\ee_1})+\sign(\xi^t_{\ii+\ee_2}))\}\;.
  \end{split}
\end{equation}
We want to check precisely which condition on $\xi^t$ ensures that
$\xi^{t+1}\in M(\xi^t)$, \ie that $(\ii,t)$ is not an error site for
$(\xi,\xx_{\Pi(\ii)})$. We know from section \ref{sec:local-coupling} that if
\begin{equation}
  \xi^t\in G
  :=
  \big\{y\in I^{U(\ii)}:|y_\ii|, |y_{\ii+\ee_1}| \text{ and } 
  |y_{\ii+e_2}| \in[\frac{d}{v},v]\big\}\;,
\end{equation}
then $\hat T\Phi_{\epsilon}^{\{t\}}(\xi^t)\in M(\xi^t)$. If we assume further
that $\hat T\Phi_{\epsilon}^{\{t\}}(\xi^t)\in G$, then the fact that
$\tilde\tau^k([\frac dv,v])\subset\tilde\tau^k([c,v])=[c,v]$ implies that also
$\hat T^2\Phi_{\epsilon}^{\{t\}}(\xi^t)\in M(\xi^t)$ and finally, by the same
argument, we have $\hat T^3\Phi_{\epsilon}^{\{t\}}(\xi^t)\in M(\xi^t)$ provided
also $\hat T^2\Phi_{\epsilon}^{\{t\}}(\xi^t)\in G$.  This can be resumed by
\begin{equation}
  \begin{split}
    E_{(\ii,t)}:=&\{\xi\in I^{U(\ii)}: (\ii,t) 
    \text{ is an error site for }(\xi,\xx_{\Pi(\ii))}\}\\
    \subseteq &
    \{\xi:\xi^t\not\in G\}
    \cup\{\xi:\hat T\Phi_\epsilon^{\{t\}}(\xi^t)\not\in G\}
    \cup \{\xi:\hat T^2\Phi_\epsilon^{\{t\}}(\xi^t)\not\in G\}\;.
  \end{split}
\end{equation}


This description will be used to estimate (uniformly in $\xx_{\Pi(\ii)}$) the
integral from (\ref{eq:space-decoupl}). Let us assume for the moment that
$t_1=\min \{t\,:\, (\ii,t)\in\tilde V_{\{\ii\}}\}\geq 1$ and define
\begin{align*}
  \tilde V^{(t_1)}_{\{\ii\}}
  &=
  \big\{(\ii,t-t_1-1)\,:\,
  t>t_1\text{ and }(\ii,t)\in \tilde V_{\{\ii\}}\big\}\\
  B_1(\xx_{\Pi(\ii)})
  &=
  \{y\in I^{U(\ii)}\,:\,(y,T_{\epsilon,\Pi(\ii)}^{t_1+1}(\xx_{\Pi(\ii)}))
  \in E(\tilde V^{(t_1)}_{\{\ii\}})\}
\end{align*}
We have then $B_0(\xx_{\Pi(\ii)})=E_{(\ii,t_1)}\cap
T^{-(t_1+1)}_{\epsilon|\xx_{\Pi(\ii)}}B_1(\xx_{\Pi(\ii)})$, which allows us to
write, for each $h_0:I^{U(\ii)}\to[0,\infty)$ of bounded variation,
\begin{equation}
  \label{eq:esterr}
  \begin{split}
    &\hspace*{-10mm}
    \int_{B_0(\xx_{\Pi(\ii)})} h_0(y)\,dy\\
    =&
    \int_{B_1(\xx_{\Pi(\ii)})} P^{t_1+1}_{\epsilon|\xx_{\pi(\ii)}}
    \Big(1_{E_{(\ii,t_1)}}(y)\,h_0(y)\Big)dy\\
    \leq&
    \int_{B_1(\xx_{\Pi(\ii)})} 
    P_{\hat T}\bigg(P^2_{\hat T}P_{\Phi_\epsilon^{\{t_1\}}} 
    \Big(1_{G^c}P_{\hat T}\big(P^2_{\hat T}P_{\Phi_\epsilon^{\{t_1-1\}}} 
    P^{t_1-1}_{\epsilon|\xx_{\pi(\ii)}} h_0\big)\Big)\bigg)(y)\,dy\\
    &+\int_{B_1(\xx_{\Pi(\ii)})} P_{\hat T}\bigg(P_{\hat T} 
    \Big(1_{G^c}P_{\hat T}\big(P_{\Phi_\epsilon^{\{t_1\}}} 
    P^{t_1}_{\epsilon|\xx_{\pi(\ii)}} h_0\big)\Big)\bigg)(y)\,dy\\
    &+\int_{B_1(\xx_{\Pi(\ii)})} P_{\hat T}\Big(1_{G^c}P_{\hat T}
    \big(P_{\hat T}P_{\Phi_\epsilon^{\{t_1\}}} 
    P^{t_1}_{\epsilon|\xx_{\pi(\ii)}} h_0\big)\Big)(y)\,dy\\
    =:&
    \int_{B_1(\xx_{\Pi(\ii)})} h_1(y)\,dy
  \end{split}
\end{equation}
and we will show at the end of this proof that
\begin{equation}
  \label{eq:rec-dens}
  \var(h_1)\leq \kappa K \var(h_0)
\end{equation}
for some constant $K$ (depending only on $\eta$, $\delta$ and $\gamma$, but
not on $k$).

We can use this to obtain by induction an estimate for the integral appearing
in (\ref{eq:space-decoupl}), starting from $h_0$, the density of the
normalized Lebesgue measure on $[\frac dv, v]^{U(\ii)}$. In the exceptional
case $t_1=0$ the first of the three integrals in the decomposition
(\ref{eq:esterr}) does not make sense, but we can use $1_{G^c}h_0=0$ instead
so that we are left with the two other integrals only.  As we only use that
$h_0\geq0$, but not that $\int h_0(y)\,dy=1$, we can reproduce the
same argument for $h_1$, with $t_2\geq 1$ (since the times in $\tilde
V_{\{\ii\}}$ have been taken at distance at least 2), obtaining $h_2$ such
that
\begin{displaymath}
  \int_{B_0(\xx_{\Pi(\ii)})} h_0(y)\,dy
  \leq
  \int_{B_1(\xx_{\Pi(\ii)})} h_1(y)\,dy
  \leq 
  \int_{B_2(\xx_{\Pi(\ii)})} h_2(y)\,dy
\end{displaymath}
and $\var(h_2)\leq\kappa K\var(h_1)\leq (\kappa K)^2 \var(h_0)$. Inductively,
we obtain densities $h_1,h_2,\dots,h_{\card\tilde V_{\{\ii\}}}$ such that
\begin{equation}
  \begin{split}
    \int_{B_0(\xx_{\Pi(\ii)})} h_0(y)\,dy
    &\leq
    \int_{I^{U(\ii)}} h_{\card\tilde V_{\{\ii\}}}(y)\,dy
    \leq \frac12\var(h_{\card\tilde V_{\{\ii\}}})\\
    &\leq (\kappa K)^{\card\tilde V_{\{\ii\}}} \frac12\var(h_0)
    \leq \Delta^{8\cdot\card\tilde V_{\{\ii\}}}    
  \end{split}
\end{equation}
if $\kappa$ is chosen small enough (by taking $k$ large enough, independently
of $K$) that $\kappa K \frac12 \var(h_0)\leq \Delta^8$. Inserted into
(\ref{eq:space-decoupl}) this yields (\ref{eq:exp-est-induct}).

To obtain (\ref{eq:rec-dens}), we notice that the three terms in the sum have
the same structure
\begin{displaymath}
  \int_{B_1(\xx_{\Pi(\ii)})} 
  P_{\hat T} P_{A_1} \Big(1_{G^c}P_{\hat T}P_{A_2} h_0\Big)(y)\,dy
\end{displaymath}
with $P_{A_1}$ and $P_{A_2}$ representing operators (depending on $t_1$ and
$\xx_{\Pi(\ii)}$) which are integral preserving and uniformly (in
$t_1$ and $\xx_{\Pi(\ii)}$) bounded in
variation norm, see (\ref{eq:LY-conditional}). As $G^c$ is the union of a
finite number of hyper-rectangles in $I^{U(\ii)}$, multiplication by $1_{G_c}$
is a bounded linear operator with respect to the norm $\var(.)$.

We denote by $K_1$ a uniform bound for the variation norm of $P_{A_1}
1_{G^c}P_{\hat T}P_{A_2}$ and by $K_2$ a uniform bound for that of $P_{A_2}$.
One can then estimate uniformly variation and integral of $P_{A_1}
1_{G^c}P_{\hat T}P_{A_2} h_0$:
\begin{equation}
  \label{eq:varcomp}
  \var(P_{A_1} 1_{G^c}P_{\hat T}P_{A_2} h_0)\leq K_1\var(h_0)
\end{equation}
and, applying (\ref{eq:JB1}) to $f=P_{A_2}h_0$,
\begin{displaymath}
  \int P_{A_1}1_{G^c} P_{\hat T}P_{A_2} h_0\,dm
  =
  \int1_{G^c} P_{\hat T}P_{A_2} h_0\,dm
  =
  \int1_{G^c}(h+\tilde h)\,dm  
\end{displaymath}
with
\begin{itemize}
\item $\var(\tilde h)\leq\kappa \var(f)\leq \kappa K_2\var(h_0)$ 
  so that $\int1_{G^c}\tilde h\,dm\leq \frac12\kappa K_2\var(h_0)$, and
\item $h=\Big(\int h_0\,dm \Big) \hat h_\alpha$ so that, according to
  (\ref{eq:kappa}),
  \begin{displaymath}
    \int 1_{G^c}h\,dm
    \leq 
    \Big(\int h_0\,dm\Big)\kappa\leq \frac 12\kappa\var(h_0)\;.    
  \end{displaymath}
\end{itemize}
Hence we obtain
\begin{equation}
  \label{eq:intcomp}
  \int P_{A_1}1_{G^c} P_{\hat T}P_{A_2} h_0\,dm
  \leq
  \frac12\kappa(1+K_2)\var(h_0)
\end{equation}
and estimates $(\ref{eq:varcomp})$ and $(\ref{eq:intcomp})$ allow us to apply
(\ref{eq:LY-1D}) to obtain
\begin{equation}
  \label{eq:vartota}
  \var(P_{\hat T} P_{A_1} 1_{G^c}P_{\hat T}P_{A_2} h_0)
  \leq 
  \kappa\big(K_1+\frac12(1+K_2)\beta\big)\var(h_0)
\end{equation}
which implies (\ref{eq:rec-dens}) with $K=3(K_1+\frac12(1+K_2)\beta)$.

\section{Phase transition versus bifurcation}
\label{sec:phasetr-bif}
\noindent\textbf{Proof of Theorem~\ref{theo:finite-dim}}\; Now 
$\Lambda=(\ZZ/L\ZZ)^2$ so that $T_\epsilon$ is a piecewise affine, piecewise
expanding self map of the $L^2$-dimensional cube $[-1,1]^\Lambda$.  The
existence of a $T_\epsilon$-invariant probability measure
$d\mu_{\epsilon,\Lambda}=h_{\epsilon,\Lambda}dm^\Lambda$ follows rather
immediately from Lemma~\ref{lem:lasota-yorke}, see \eg Theorem~3.1 in
\cite{KeLi-2005a}. 
Indeed, more is true: As the Perron-Frobenius operator of $T_\epsilon$ is
quasicompact, there is some iterate $T_\epsilon^r$ which is exact on
each of its ergodic components. Therefore the uniqueness and mixing of
$\mu_{\epsilon,\Lambda}$ can be proved
following
the folklore type strategy of \cite{Keller-1997}. Namely, we will show:
\begin{equation}
  \label{eq:B}
  \parbox{10cm}{If $B\subseteq I^\Lambda$ is a measurable 
    $T_\epsilon^r$-invariant set,\\ then
    $B=I^\Lambda$ modulo Lebesgue measure $0$.}
\end{equation}
This proves uniqueness of the invariant density and mixing of
$(T_\epsilon,\mu_{\epsilon,\Lambda})$ and, in view of the spectral gap of the
Perron-Frobenius operator of $T_\epsilon$, also the exponential decay of
correlations follows at once.  In the sequel we will suppress to write
``modulo Lebesgue measure $0$'', but all set inclusions will be understood in
this way.

Denote by $\J$ the partition of $I=[-1,1]$ into maximal monotonicity intervals
of $\hat\tau$. Then $\J^\Lambda$ is a partition of $I^\Lambda$ into
hyperrectangles $R$ each of which is mapped by $\hat T$ homeomorphically onto
its image. As some maximal monotone branches of $\hat\tau$ are only piecewise
linear, also the $\hat T|_R$ are only piecewise linear in general.
Nevertheless, $\hat T(R)$ is a hyperrectangle for each $R\in\J^\Lambda$.

Now suppose $B$ is as in \eqref{eq:B}. A telescoping argument just as in
\cite[proposition 5]{Keller-1997} shows that $B$ contains some cube $Q$ of
side length $\ell(Q)$.  Elementary geometrical arguments (see \cite[Lemma
8a]{Keller-1997}) show that $\Phi_\epsilon Q$ contains a cube $Q'$ of side
length at least $\frac\epsilon{1-\epsilon}\ell(Q)\geq\frac13\ell(Q)$ (observe
that $\epsilon\in[0,\frac14]$). This cube $Q'$ is partitioned by the partition
$\J^\Lambda$ into hyperrectangles $Q'\cap R$. As long as $Q'$ is cut in no
more than two pieces in each coordinate direction, at least one of the $Q'\cap
R$ contains a cube $Q''$ of side length at least $\frac12\ell(Q')$.  As we
chose the parameter $k$ in Lemma~\ref{lem:lasota-yorke} so large that
$|\hat\tau'|\geq4$, this guarantees that $\ell(\hat TQ'')\geq2\ell(Q')$. In
this way we conclude that $\hat T^3 Q'$ contains a subcube $Q'''$ of side
length at least $8\ell(Q')\geq\frac83\ell(Q)$, provided that by none of the
three iterated applications of $\hat T$ any side of the cube is cut into more
than two pieces.

If, on the other hand, at any of the three stages some sides are cut into more
than two pieces, then one of these pieces is a maximal monotonicity interval
of $\hat\tau$, and as we saw in the proof of Lema~\ref{lem:single-site-1}, its
image has length at least $\frac{u-c}v$.  Hence, it is still guaranteed that
$T_\epsilon Q\supseteq\hat T^3Q'$ contains a cube $Q'''$ with
$\ell(Q''')\geq\min\{\frac83\ell(Q),\frac{u-c}v\}$.  Continuing in this way,
$T_\epsilon^{nr}Q$ will contain a cube of side length at least $\frac{u-c}{v}$
for some $n>0$. As $Q\subset B=T_\epsilon^rB$, we conclude that $B$ contains a
cube of side length at least $\frac{u-c}{v}$, and so $\Phi_\epsilon B$
contains a cube, call it $\hat Q$, with $\ell(\hat Q)\geq\frac{u-c}{3v}$

Now, if $k$ is chosen large enough, then each monotonicity interval of
$\hat\tau=\frac1v\tilde\tau^k$ has length less than $\frac{u-c}{6v}$.
Therefore each side of $\hat Q$ contains at least one maximal monotonicity
interval of $\hat\tau$, and it follows from Lemma~\ref{lem:single-site-1} that
$T_\epsilon B=\hat T^3\Phi_\epsilon B\supseteq\hat T^3\hat Q=I^\Lambda$. Hence
$B=T_\epsilon^rB=T_\epsilon^{r-1}T_\epsilon
B=T_\epsilon^{r-1}I^\Lambda=I^\Lambda$, \ie \eqref{eq:B}. This proves at the
same time that the invariant measure is indeed equivalent to Lebesgue measure.

\section{The phase transition for smooth expanding circle maps}
\label{sec:smooth}

In this section we describe how to modify our construction such that $\tau$
becomes a smooth expanding circle map. We follow essentially the strategy from
\cite{KeRu-2004} where a piecewise linear map, whose Perron-Frobenius operator
had particular spectral properties, was approximated by an analytic map in such
a way that the spectal properties were essentially conserved. Our situation is
different from the one there in so far as the map $\tau$ has non-surjective
monotone branches, so a little modification of the strategy from
\cite{KeRu-2004} is necessary.
\\[3mm]
\textbf{The construction of the modification $\bar\tau$ of $\hat\tau$}\\
Denote by $\acute\tau$ the modification of $\tilde\tau$ with only increasing
linear branches as discussed in
Remark~\ref{rem:full-branches}(\ref{item:increasing}). So
Lemma \ref{lem:single-site-2} holds as well for the map
$\check\tau:=\frac1v\acute\tau^k$. Let $-1=a_0<a_1<\dots<a_p=1$ be the partition
of $[-1,1]$ into maximal monotonicity intervals of $\check\tau$. We are going
to define an increasing, piecewise linear, continuous function
$\zeta:\RR\to\RR$ in terms of the inverse branches of $\check\tau$: for
$k=0,\dots,p-1$ and $x\in(2k-1,2k+1)$ let
\begin{equation}
  \label{eq:zeta-def1}
  \zeta(x)=
  \begin{cases}
    a_{k}&\text{ if }x-2k\leq\check\tau(a_{k}^+)\\
    \check\tau|_{(a_{k},a_{k+1})}^{-1}(x-2k)&\text{ if
    }\check\tau(a_{k}^+)<x-2k<\check\tau(a_{k+1}^-)\\
    a_{k+1}&\text{ if }x-2k\geq\check\tau(a_{k+1}^-)\;.
  \end{cases}
\end{equation}
Obviously $\zeta$ extends continuously to a (not strictly!) increasing map
from $[-1,2p-1]$ onto $[-1,1]$. As $\zeta(-1)=a_0=-1$ and $\zeta(2p-1)=a_p=1$,
it further extends to a continuous increasing map $\zeta:\RR\to\RR$ by
\begin{equation}
  \label{eq:zeta-def2}
  \zeta(x+2p)=\zeta(x)+2\;.
\end{equation}
Let $\vf_\sigma:\RR\to[0,\infty)$ be the Gaussian density with mean $0$ and
variance $\sigma^2$. We use it as a convolution kernel for defining
\begin{equation}
  \label{eq:zeta-sigma-def}
  \zeta_\sigma(x):=\zeta*\vf_\sigma(x)-\zeta*\vf_\sigma(-1)-1\;.
\end{equation}
Obviously, $\zeta_\sigma(-1)=-1$, $\zeta_\sigma(x+2p)=\zeta_\sigma(x)+2$,
$\zeta_\sigma'(x+2p)=\zeta_\sigma'(x)$, and as $\zeta$ is continuous and
piecewise linear, $\|\zeta_\sigma-\zeta\|_\infty\leq\const\cdot\sigma$ and
$\lim_{\sigma\to0}\zeta_\sigma'(x)=\zeta'(x)$ at all points $x$ where $\zeta$
is differentiable. Furthermore, $\inf_x\zeta_\sigma'(x)>0$ for each
$\sigma>0$. Hence $\zeta_\sigma^{-1}$ projects to a $p$-fold covering circle
map $\check\tau_\sigma$ of $\RR/(2\ZZ+1)$ onto itself. In this way we are
nearly in the same situation as in \cite{KeRu-2004} except that
$\inf_{\sigma>0}\inf_x\zeta_\sigma'(x)=\inf_x\zeta'(x)=0$.\footnote{Our maps
  $\zeta$ and $\zeta_\sigma$ play the roles of the maps $\tau$ and
  $\tau_\delta$ from \cite[section 3]{KeRu-2004}. Note also that we are
  dealing with maps of $[-1,1]$ whereas \cite{KeRu-2004} considers maps of
  $[0,1]$ or $\RR/\ZZ$. We also take this opportunity to correct two misprints
  in \cite{KeRu-2004}. The first one is obvious: above eq.~(2) one has
  $\dot\tau_\delta(x+p)$ instead of $\dot\tau_\delta(x+1)$. The second one
  concerns eq.~(3), which is not correct as stated. It should be
  $\|2g_{\delta,M}\|_\infty^{1/M}\leq2^{1/M}\cdot\|g\|_\infty$, and in the
  line thereafter one must restrict to $\|g\|_\infty<\kappa<1$ instead of
  $\vartheta<\kappa<1$. This does not affect the main result of
  \cite{KeRu-2004} because, for the particular map studied there,
  $\vartheta=\|g\|_\infty$. Also for the purposes of the present paper this
  weaker form of eq.~(3) is suficient.}  The Perron-Frobenius operator of
$\check\tau_\sigma$ can be conveniently expressed in terms of $\zeta_\sigma$:
for $x\in[-1,1)$,
\begin{displaymath}
  P_{\check\tau_\sigma}f(x)
  =
  \sum_{k=0}^{p-1}f(\zeta_\sigma(x+2k))\cdot\zeta_\sigma'(x+2k)\;.
\end{displaymath}
\\[3mm]
\textbf{Proof of Theorem~\ref{theo:main} for the map $\bar\tau^3$}\\
Obviously, $\lim_{\sigma\to0}P_{\check\tau_\sigma}f(x)=P_{\check\tau}f(x)$ at
all points $x$ where $f$ is continuous and $\zeta$ is differentiable. Hence,
although $\check\tau$ has many non-surjective branches whereas
$\check\tau_\sigma$ has only full branches, $P_{\check\tau_\sigma}$ is
``close'' to $P_{\check\tau}$. Below we state this more precisely: we will
show that, for small enough $\sigma>0$, the Perron-Frobenius operator of
$\bar\tau:=\check\tau_\sigma$ is only a small modification of that of
$\check\tau$ in the sense discussed in
Remark~\ref{rem:full-branches}(\ref{item:smooth}), namely: there are constants
$\sigma_0,\kappa,F,C_0>0$ such that, for all $f:[-1,1]\to\RR$ of bounded
variation, for all $\sigma\in(0,\sigma_0)$ and for all $n\in\NN$,
\begin{gather}
  \label{eq:KR1}
  \var(P_{\check\tau_\sigma}^nf)
  \leq
  F\cdot\left(\kappa^n\var(f)+\int|f|\,dm\right)\;,\\
  \label{eq:KR2}
  \int|P_{\check\tau_\sigma}f-P_{\check\tau}f|\,dm
  \leq
  C_0\cdot\sigma\cdot\var(f)\;.
\end{gather}
Hence Theorem~\ref{theo:main} also holds for the single-site map
$\tau=\bar\tau^3$ (see Remark~\ref{rem:full-branches}(\ref{item:smooth})).

The uniform Lasota-Yorke type estimate \eqref{eq:KR1} is proved exactly as in
\cite{KeRu-2004} (which is based in turn on Rychlik's approach
\cite{Rychlik-1983}).\footnote{As in our case $\check\tau'\geq4$, one can make
  the following simplifying choices in \cite[section 2]{KeRu-2004}: $N=M=1$,
  $\epsilon=\frac14$, $\lambda_1=\kappa=\frac34$ and
  $D_1=2/\min_{A\in\alpha}m(A)$ if $\sigma>0$ is small enough, where the
  partition $\alpha=\alpha_1$ is obtained from the partition into maximal monotonicity
  intervals of $\check\tau$ by cutting each of these intervals into two subintervals with equal length, and finally $F=4D_1$.}  The proof of
\eqref{eq:KR2} is a bit different from that in \cite{KeRu-2004}, because,
other than in that reference,
$\inf_{\sigma>0}\inf_x\zeta_\sigma'(x)=\inf_x\zeta'(x)=0$. Denote by $\chi$  the
indicator function of an interval $[-1,x]$. As in
\cite[Proposition~2]{KeRu-2004} it suffices to show that
$\int|P_{\check\tau_\sigma}\chi-P_{\check\tau}\chi|\,dm\leq C_0\,\sigma$.  But
$\lim_{s\to0}P_{\check\tau_s}f(x)=P_{\check\tau}f(x)$ at all but at most
countably many $x$, so
$\lim_{s\to0}\int|P_{\check\tau_s}\chi-P_{\check\tau}\chi|\,dm=0$. Therefore
the following estimate, which is uniform in $s\in(0,\sigma)$, completes the
proof of \eqref{eq:KR2}:
\begin{equation}
  \begin{split}
    &\int|P_{\check\tau_\sigma}\chi-P_{\check\tau_s}\chi|\,dm\\
    \leq&
    \int_{-1}^1\sum_{k=0}^{p-1}\left|\chi(\zeta_\sigma(y+2k))\cdot\zeta_\sigma'(y+2k)-\chi(\zeta_s(y+2k))\cdot\zeta_s'(y+2k)\right|\,dy\\
    =&
    \int_{-1}^{2p-1}\left|\chi(\zeta_\sigma(y))\cdot\zeta_\sigma'(y)-\chi(\zeta_s(y))\cdot\zeta_s'(y)\right|\,dy\\
    \leq&
    2\|\zeta_\sigma-\zeta_s\|_\infty
    +
    \int_{-1}^{2p-1}\left|\zeta_\sigma'(y)-\zeta_s'(y)\right|\,dy
  \end{split}
\end{equation}
and both integrals are bounded by $\const\cdot\sigma$ - the first one, as
remarked above, because $\zeta$ is Lipschitz continuous, and the second one
because $\zeta'$ is of bounded variation.
\\[3mm]
\textbf{Proof of Theorem~\ref{theo:finite-dim} for the map $\bar\tau^3$}\\
Let $\bar\tau=\check\tau_\sigma$ with a sufficiently small but fixed
$\sigma>0$ be as before.  Although the branches of $\bar\tau$ are not linear,
the uniform expansion in conjunction with the bounded second derivative of
$\bar\tau$ yields uniform distortion control on the branches of all iterates
$\bar\tau^n$, \cf \cite[Lemma 5.1.18]{KaHa}. This allows the same telescoping
argument as in section~\ref{sec:phasetr-bif} leading to a cube $Q$ as in that
proof. The rest of the proof goes through exactly as in that section. Indeed,
the argument becomes even much simpler, because now all monotone branches of
the single-site map are surjective.

\end{document}